\documentclass[12pt]{article}

\usepackage{amsmath,amssymb,amsfonts,theorem,makeidx,latexsym,epsfig,,subfigure}

\newtheorem{defn}{Definition}[section]

\newtheorem{lemma}[defn]{Lemma}

\newtheorem{proposition}[defn]{Proposition}

{\theorembodyfont{\rmfamily}

\newtheorem{ex}[defn]{Example}}

\newtheorem{thm}[defn]{Theorem}

\newtheorem{prop}[defn]{Proposition}

\newtheorem{cor}[defn]{Corollary}

\numberwithin{equation}{section}

\newcommand{\h}{{\cal H}}

\newcommand{\ltr}{ L^2(\mathbb R) }

\newcommand{\mn}{\mathbb N}

\newcommand{\mr}{\mathbb R}

\newcommand{\mz}{\mathbb Z}

\newcommand{\mc}{\mathbb C}

\newcommand{\mt}{ E_{mb}T_{na}g}

\def\bp{{\noindent\bf Proof. \ }}

\def\ep{\hfill$\square$\par\bigskip}

\def\bqs{\begin{equation}}

\def\eqs{\tag*{$\square$}\end{equation}\par\bigskip}

\def\la{\langle}

\def\ra{\rangle}

\def\ga{\gamma}

\def\bop{\begin{op}\rm}

\def\eop{\end{op}}

\def\bee{\begin{eqnarray}}

\def\ene{\end{eqnarray}}

\def\bes{\begin{eqnarray*}}

\def\ens{\end{eqnarray*}}

\def\bei{\begin{itemize}}

\def\eni{\end{itemize}}

\def\bt{\begin{thm}}

\def\et{\end{thm}}

\def\bc{\begin{cor}}

\def\ec{\end{cor}}

\def\bpr{\begin{prop}}

\def\epr{\end{prop}}

\def\bl{\begin{lemma}}

\def\el{\end{lemma}}

\def\bd{\begin{defn}}

\def\ed{\end{defn}}

\def\bex{\begin{ex}}

\def\enx{\end{ex}}

\def\bfi{\begin{fig}}

\def\efi{\end{fig}}

\def\inr{\int_{-\infty}^\infty}

\newcommand{\nft}{ || f||^2}

\def\mt{{\mathbb T}}

\def\fti{\{f_i\}_{i\in I}}
\def\gti{\{g_i\}_{i\in I}}

\def\sui{\sum_{i\in I}}

\title{On extensions of wavelet systems to dual pairs of frames}

\date{\today}

\author{Ole Christensen\thanks{Department of Applied Mathematics and Computer Science,
Technical University of Denmark,
Building 303, 2800 Lyngby, Denmark (ochr@dtu.dk)
},
Hong Oh Kim\thanks{Department of Mathematical
Sciences, KAIST, 373-1, Guseong-dong, Yuseong-gu, Daejeon, 305-701,
Republic of Korea (kimhong@kaist.edu)},
Rae Young Kim\thanks{Department of Mathematics,
Yeungnam University,
214-1, Dae-dong, Gyeongsan-si, Gyeongsangbuk-do, 712-749,
Republic of Korea (rykim@ynu.ac.kr)
}}

\begin{document}

\maketitle

\begin{abstract}
It is an open problem whether any pair of Bessel sequences with
wavelet structure can be extended to a pair of dual frames by
adding a pair of singly generated wavelet systems. We consider the
particular case where the given wavelet systems are generated by
the multiscale setup with trigonometric masks and provide a
positive answer under extra assumptions. We also identify a number
of conditions that are necessary for the extension to dual
multiscale wavelet frames with any number of generators, and show
that they imply that an extension with two pairs of wavelet
systems is possible. Along the way we provide examples showing
that extensions to dual frame pairs are attractive because they
often allow better properties than the more popular extensions to
tight frames.
\end{abstract}

\begin{minipage}{150mm}
\end{minipage}


\section{Introduction}

 Extension problems exist in a large variety in the frame
literature. In its basic version the question is whether a given
sequence of elements in a Hilbert space can be extended to a frame
with prescribed properties. For example, it is natural to ask for
extensions such that the resulting frame is computationally
convenient, e.g., a tight frame or a frame for which a dual frame
can be found easily. A natural generalization of this idea is to
start with {\it two} sequences and ask for extension of these
sequences to dual frame pairs.

It is known that any pair of Bessel sequences in a separable
Hilbert space can be extended to a pair of dual frames by adding
appropriate collections of vectors. But if we require the added
sequences to have a special structure or to satisfy certain
constraints, many open problems appear. A key question is whether
any given pair of Bessel sequences with wavelet structure can be
extended to a pair of dual frames by adding  a pair of wavelet
systems, each with a single generator (see the paper \cite{DH} by
D. Han for a discussion and conjecture about the corresponding
tight case).

In this paper we will analyze this problem under the extra
assumption that the given wavelet systems are generated by the
MRA-setup as considered, e.g., in \cite{RS97-2, CHS, DRoSh3}. Given two scaling functions, we will
consider the associated wavelet systems generated by letting the
masks be trigonometric polynomials, and ask for extensions to dual
pairs of frames by adding wavelet systems of the same type.  We
will first identify a condition on the refinement masks that is
necessary for this extension to be possible, and then show that
this condition is also sufficient for the possibility to extend to
dual pairs using two generators. A stronger condition
characterizes the possibility to extend the given wavelet systems
to a pair of dual wavelet frames by adding a pair of wavelet
systems, each with a single generator.

Note that Daubechies and B. Han  already in \cite{DH04}   showed
that for any given pair of scaling functions one can construct
dual pairs of wavelet frames, each with two generators. Our setup
is different from the one in \cite{DH04}: in our extension of
Bessel systems to dual pairs with two generators, we consider  one
of the pairs of Bessel sequences to be given, i.e., we only have
freedom with respect to one pair of wavelet systems.

The rest of  this introduction gives a short introduction to the
key ingredients of the paper and connects our results to the
literature. In Section \ref{1812a} we present a few preliminary
results on Bessel sequences and the mixed extension principle. The
results about extensions with singly generated systems are stated
in Section \ref{1812b}, while the extension by systems with two
generators is treated in Section \ref{1812c}.

A sequence $\fti$ in a separable Hilbert space $\h$ is called a
{\it  frame} if there exist constants $A,B>0$ such that \bes A\,
\nft \le \sui |\la f, f_i \ra|^2 \le B\,\nft, \ \forall f\in
\h.\ens A frame is {\it tight} if we can choose $A=B;$ and if at
least the upper frame condition is satisfied, $\fti$ is called a
{\it Bessel sequence.} For any frame $\fti,$ there exist at least
one {\it dual frame,} i.e., a frame $\gti$ such that $f= \sui \la
f, g_k\ra f_k, \ \forall f\in \h.$ A tight frame with $A=B=1$
leads to an expansion of arbitrary elements $f\in \h$ of exactly
the same type as we know for orthonormal bases, i.e., $f= \sui \la
f, f_k\ra f_k, \ \forall f\in \h.$ For more information on frames
we refer to the books \cite{Dau, CB}.

In this paper we will exclusively consider systems of functions in
$\ltr$ with {\it wavelet structure,} that is, collections of
functions of the type $\{2^{j/2} \psi(2^jx-k)\}_{j,k\in \mz}$ for
a fixed function $\psi.$ Considering the operators on $\ltr$ given
by $T_kf(x):=f(x-k)$ and $Df(x):= 2^{1/2}f(2x),$ the wavelet
system can be written as $\{D^jT_k\psi\}_{j,k\in \mz}.$

Let $\mt$ denote  the unit circle which will be identified with
$[-1/2,1/2].$  Also, for $f\in L^1(\mr)\cap L^2(\mr)$ we denote
the Fourier transform by ${\mathcal F}f(\ga)= \hat{f}(\ga)= \inr
f(x)e^{-2\pi i x  \ga}dx.$ As usual, the Fourier transform is
extended to a unitary operator on $\ltr$.

In the entire paper we will use the following setup that appeared,
e.g., in \cite{DRoSh3}, except that we restrict our attention to
trigonometric masks.

\vspace{.1in}\noindent{\bf General setup:} Consider a function
$\varphi\in \ltr$ such that

\bei \item[(i)] $\widehat \varphi$ is continuous at the origin and
$\hat\varphi(0)=1;$ \item[(ii)] There exists a $1$-periodic
trigonometric polynomial  $m_0$ (called a     {\it refinement
mask})
    such that \bee \label{1312a} \widehat\varphi(2\ga)=m_0(\ga)\widehat\varphi(\ga), \quad\text{a.e. } \ga\in\mr.\ene \eni
Given $1$-periodic trigonometric polynomials  $m_1,m_2, \dots,
m_n,$ consider the functions $\psi_\ell\in \ltr$ defined by \bee
\label{r-1}\widehat{\psi_\ell}(2\ga)=m_\ell(\ga)\widehat\varphi(\ga),
\ \ell=1, \dots, n.\ene

\vspace{.1in}\noindent Note that the technical condition
$$\sum_{k\in \mz}\left|  \widehat \varphi (\ga+k)  \right|^2 \leq K <\infty,$$
which is used in \cite{DRoSh3}, automatically is satisfied in our
setting. In fact, it is well known that $\varphi$ has compact
support whenever the scaling equation \eqref{1312a} holds for a
trigonometric polynomial $m_0.$ Thus,
$$\sum_{k\in \mz}\left| \widehat \varphi (\ga+k)  \right|^2
=\sum_{k\in \mz}  \langle\varphi, \varphi(\cdot -k) \rangle
e^{-2\pi i k \ga} $$ is a trigonometric polynomial. \vspace{.1in}

We will base the analysis on the {\it mixed  extension principle
}(MEP) by Ron and Shen \cite{RS97-2} which is formulated in terms
of the $(n+1)\times 2$ matrix-valued functions $M$ and $\tilde M$
defined by
\begin{equation} \label{m-12}
M(\ga)=\begin{pmatrix}
  m_0\left(\ga \right) & m_0\left(\ga+\frac{1}{2}\right) \\
  m_1\left(\ga \right) & m_1\left(\ga+\frac{1}{2}\right) \\
  \vdots & \vdots  \\
  m_n\left(\ga \right) & m_n\left(\ga+\frac{1}{2}\right)
\end{pmatrix},
\quad \widetilde M(\ga)=\begin{pmatrix}
\widetilde  m_0\left(\ga \right) &\widetilde m_0\left(\ga+\frac{1}{2}\right) \\
\widetilde  m_1\left(\ga \right) &\widetilde m_1\left(\ga+\frac{1}{2}\right) \\
  \vdots & \vdots  \\
\widetilde  m_n\left(\ga \right) &\widetilde
m_n\left(\ga+\frac{1}{2}\right)
\end{pmatrix}.
\end{equation}

Formulated for trigonometric masks, the MEP reads as follows:

\begin{proposition}\label{MEP}
Assume that $\varphi, \widetilde \varphi \in L^2(\mr)$ satisfy the
conditions in the general setup, with associated masks  $m_0,
\widetilde m_0.$ For each $\ell=1,\cdots,n$, let $m_\ell,
\widetilde m_\ell$ be trigonometric polynomials and define
$\psi_\ell, \widetilde \psi_\ell\in L^2(\mr)$ by \eqref{r-1}.
Assume that $\{D^jT_k\psi_i\}_{i=1,\cdots,n; j, k \in \mz}$ and
$\{D^jT_k\widetilde \psi_i\}_{l=1,\cdots,n; j, k \in \mz}$ are
Bessel sequences. If the corresponding matrix-valued functions $M$
and $\widetilde M$ satisfy
\begin{equation} \label{m-11}
\widetilde M(\ga)^* M(\ga)=I, \  \ \ga\in \mt,
\end{equation}
then $\{D^jT_k\psi_i\}_{i=1,\cdots,n; j, k \in \mz}$ and
$\{D^jT_k\widetilde \psi_i\}_{l=1,\cdots,n; j, k \in \mz}$ form
dual frames for $L^2(\mr)$.
\end{proposition}

The MEP was later extended to the {\it mixed oblique extension
principle} in  \cite{CHS, DRoSh3}; these papers also
contain several explicit examples. Other papers about the MEP
include \cite{BGN, E, HS}.

Most of the concrete wavelet frame constructions in the literature
are obtained  via the related unitary extension principle
\cite{RS97-2} and its variants, which lead to tight frame
constructions (see, e.g., the papers \cite{RoSh4, P, HM3, CPSS, CKK-CMB}, just to mention a few out
of many). But it is already noted in, e.g., \cite{DRoSh3, CHS} that the extra flexibility in the MEP frequently leads
to more attractive constructions, a claim that is also supported
by some of the results in the current paper. For example, we
consider a case  where the extension of a wavelet system to a
tight frame introduces a wavelet generator without compact
support, while the extension to a dual pair of frames is possible
with compactly supported generators.

Note that the  analysis in the current paper is complementary to
the one in \cite{CKK}. In \cite{CKK} we formulated the general
question whether any pair of Bessel sequences
$\{D^jT_k\psi_1\}_{j, k \in \mz}$ and $\{D^jT_k \widetilde
\psi_1\}_{j, k \in \mz}$ can be extended to a pair of dual frames
by adding a wavelet system to each of the given Bessel sequences.
A sufficient condition for a positive answer turned out to be that
$\widehat{\psi_1}$ is compactly supported on $[-1,1].$ In
contrast, the extension principle applied in the current paper
involves functions that are compactly supported in time.

Finally, for the sake of the non-specialist, we note that it is
known that the dual frames of a wavelet frame not necessarily have
wavelet structure: there are cases where no dual wavelet frame
exist at all (see, e.g., the books by Chui \cite{Chui} and
Daubechies \cite{Dau}), and there are cases where some duals have
wavelet structure and some do not (see the paper by Bownik and
Weber \cite{BW}).  This issue is one of the key motivations behind
the various extension principles, which construct, simultaneously,
a frame and a dual with wavelet structure.

\section{Preliminaries on Bessel sequences and the  MEP} \label{1812a}

In the entire paper we assume that we have given  trigonometric
polynomials $m_0,m_1, \widetilde m_0, \widetilde m_1 $ as
described in the general setup. We will search for  trigonometric
polynomials $m_2, \cdots, m_n$ and $ \widetilde m_2, \cdots,
\widetilde m_n$ for  the cases $n=2$ and $n=3$ satisfying the
condition \eqref{m-11}, $i.e.$,
\begin{equation} \label{1212a}
\widetilde M(\ga)^* M(\ga)=I, \  \ \ga\in \mt.
\end{equation}
Note that \eqref{1212a} is equivalent to the two conditions
\begin{eqnarray}
   && \sum_{\ell=0}^n \overline{m_\ell(\ga)}\widetilde m_\ell(\ga)=1, \  \ \ga\in \mt,\label{m-2}\\
&& \sum_{\ell=0}^n \overline{m_\ell(\ga)}\widetilde
m_\ell(\ga+1/2)=0,\  \ \ga\in \mt. \label{m-3}
\end{eqnarray}

In Proposition \ref{MEP} we need that
$\{D^jT_k\psi_i\}_{i=1,\cdots,n; j, k \in \mz}$ and
$\{D^jT_k\widetilde \psi_i\}_{l=1,\cdots,n; j, k \in \mz}$ are
Bessel sequences. Since these systems are finite union of wavelet
systems, it is sufficient that each of these form a Bessel
sequence. The following lemma provides necessary and sufficient
conditions for this:

\bl \label{m-1} Let $\varphi \in L^2(\mr)$ be a scaling function,
with a refinement mask $m_0$  satisfying the conditions in the
general setup. Let $m$ be a trigonometric polynomial and define
$\psi$ by
 \bee
\label{s-1}\widehat{\psi}(2\ga)=m(\ga)\widehat\varphi(\ga), \
\ell=1, \dots, n.\ene
Then the following are equivalent:
\begin{itemize}
    \item[{\rm (a)}] $\{D^jT_k\psi\}_{ j, k \in \mz}$ is a Bessel sequence;
    \item[{\rm (b)}] $m(0)=0.$
\end{itemize}
\el

\bp (a) $ \Rightarrow$ (b) : Assume that $\{D^jT_k\psi\}_{ j, k
\in \mz}$ is a Bessel sequence with bound $B$. Using \cite[Theorem
1]{CS93} (or see Proposition 11.2.2 in \cite{CB}), we have
$\sum_{j\in \mz}\left| \widehat \psi (2^j \ga)  \right|^2 \leq B.$
By the conditions in the general setup, $\widehat \psi$ is
continuous at the origin, so it follows that $\widehat \psi(0)=0.$
Hence, via \eqref{s-1} and the assumption $\widehat \varphi(0)=1$
we conclude  that $m(0)=0$.

\noindent (b) $ \Rightarrow$ (a) : See \cite[Lemma 2.1]{DH04}. \ep

\bex Let $B_N, \ N\in \mn,$ denote the $N$th order B-spline,
defined recursively by \bes B_1:= \chi_{[0,1]}, \ B_{N+1}:= B_1*
B_N.\ens It is well known that $\varphi:=B_N$ satisfies the
conditions in the general setup with $m_0(\ga)= \left(e^{-\pi i
\ga} \cos(\pi \ga)\right)^N.$ Thus Lemma \ref{m-1} shows that
$\{D^jT_kB_N\}_{j,k\in \mz}$ is not a Bessel sequence in $\ltr.$
On the other hand, a finite linear combination \bes \psi_1(x)=
\sum_{k=-K}^K c_k B_N(2x-k)\ens generates a Bessel sequences
$\{D^jT_k\psi_1\}_{j,k\in \mz}$ if and only if $\sum_{k=-K}^K
c_k=0.$ \ep \enx

The following example demonstrates that the matrix condition in
Theorem \ref{MEP} is not sufficient for duality, i.e., the
assumption of the sequences $\{D^jT_k\psi_i\}_{i=1,\cdots,n; j, k
\in \mz}$ and $\{D^jT_k\widetilde \psi_i\}_{l=1,\cdots,n; j, k \in
\mz}$ being Bessel sequences is essential:
\begin{ex}
Consider the scaling functions $\varphi(x)=\widetilde
\varphi(x)=B_2(x)$ and
\begin{eqnarray*}
   &&\psi_1(x)=\widetilde \psi_1(x)=-\frac{1}{2}B_2(2x) +B_2(2x-1) -\frac{1}{2}B_2(2x-2),\\
&& \psi_2(x)=-\frac{1}{4}B_2(2x+2) +\frac{1}{2}B_2(2x)
-\frac{1}{4}B_2(2x-2), \ \ \widetilde \psi_2(x)=2B_2(2x),
\end{eqnarray*}
with the associated masks
\begin{eqnarray*}
  && m_0(\ga)= \widetilde m_0(\ga)=  e^{-2\pi i \ga} \cos^2(\pi \ga), \\
 && m_1(\ga)=  \widetilde m_1(\ga)=  e^{-2\pi i \ga} \sin^2(\pi \ga), \ \
m_2(\ga)=  2\cos^2(\pi \ga) \sin^2(\pi \ga),\ \  \widetilde
m_2(\ga)= 1.
\end{eqnarray*}
Then the MEP-condition \eqref{m-11} with $n=2$  in Proposition
\ref{MEP} is satisfied. But  by Lemma \ref{m-1} we know that
$\{D^jT_k \widetilde \psi_2\}_{j, k \in \mz}$ is not a Bessel
sequence. Hence $\{D^jT_k\psi_i\}_{i=1,2; j, k \in \mz}$ and
$\{D^jT_k\widetilde \psi_i\}_{l=1,2; j, k \in \mz}$ does not form
dual frames for $L^2(\mr)$. \ep
\end{ex}

In Lemma \ref{m-5} we will state three necessary conditions for
the existence of MEP-type wavelet systems
$\{D^jT_k\psi_i\}_{i=1,\cdots,n; j, k \in \mz}$ and
$\{D^jT_k\widetilde \psi_i\}_{l=1,\cdots,n; j, k \in \mz}$. We
first need the following factorizations:

\bl \label{m-4} \

\begin{itemize}
    \item[{\rm(i)}] Let $f$ be a $1$-periodic trigonometric polynomial with
$f(0)=0.$ Then  $f(\ga)=e^{-\pi i \ga}\sin(\pi \ga)
\Lambda_1(\ga)$
                    for a $1$-periodic trigonometric polynomial $\Lambda_1$;
    \item[{\rm(ii)}] Let $g$ be a $1$-periodic trigonometric polynomial with
$g(1/2)=0.$ Then $g(\ga)=e^{-\pi i \ga}\cos(\pi \ga)
\Lambda_2(\ga)$
                    for a $1$-periodic trigonometric polynomial $\Lambda_2$.
\end{itemize}
\el \bp For the proof of (i), write $f$  as $f(\ga)=\sum_{k \in
\mz} c_k e^{-2\pi i k \ga}.$ Since $f(0)=\sum_{k \in \mz} c_k=0$,
we have $f(\ga)= \sum_{k \neq 0} c_k(e^{-2\pi i k \ga} -1).$
Define $f_+$ and $f_-$ by
$$f_+(\ga):=\sum_{k\in \mn } c_k (e^{-2\pi ik \ga}-1), \ \ \ \
f_-(\ga):=\sum_{k\in \mn } c_{-k} (e^{2\pi k i \ga}-1).$$ Then we
see that
\begin{eqnarray*}
   f_+(\ga)=\sum_{k\in\mn}c_k(e^{-2\pi i \ga} -1)\sum_{\ell=0}^{k-1}e^{-2\pi i \ell \ga}
   &=& e^{-\pi i \ga}\sin(\pi \ga)\left(-2i \sum_{k\in \mn} c_k \sum_{\ell=0}^{k-1} e^{-2\pi i \ell \ga}  \right)\\
         &=:& e^{-\pi i \ga}\sin(\pi \ga) \Lambda_+(\ga).
\end{eqnarray*}
Similarly,
\begin{eqnarray*}
f_-(\ga)&=&e^{-\pi i \ga}\sin(\pi \ga) \left(2i \sum_{k\in \mn}
c_{-k} \sum_{\ell=1}^{k} e^{2\pi i \ell \ga}  \right) =: e^{-\pi i
\ga}\sin(\pi \ga) \Lambda_-(\ga).
\end{eqnarray*}
Then we have $f(\ga)=f_+(\ga)+f_-(\ga)=e^{-\pi i \ga}\sin(\pi \ga)
\Lambda_1(\ga),$ where
$\Lambda_1(\ga):=\Lambda_+(\ga)+\Lambda_-(\ga)$ is a $1$-periodic
trigonometric polynomial. This proves (i).

For the proof of (ii), let $\tilde g(\ga):=g(\ga+1/2).$ Since
$\tilde g(0)=0$, there exists a $1$-periodic trigonometric
polynomial $\Lambda$ such that $\tilde g(\ga)=e^{-\pi i
\ga}\sin(\pi \ga) \Lambda(\ga).$ Then we have $g(\ga)= \tilde
g(\ga-1/2)
 = e^{-\pi i \ga} \cos (\pi \ga) \Lambda_2(\ga),$
where $\Lambda_2(\ga):=-i \Lambda(\ga-1/2).$ This proves (ii). \ep

\bl \label{m-5} Under the hypothesis of Proposition \ref{MEP}, the
following hold:
\begin{itemize}
    \item[{\rm (a)}] $m_\ell(0)=\widetilde m_\ell(0)=0$, $\ell=1,2,\cdots,n$;
   \item[{\rm (b)}] $m_0(1/2)=\widetilde m_0(1/2)=0$;
    \item[{\rm (c)}]
$ 1-\overline{m_0(\ga)} \widetilde m_0(\ga)= \sin^2\left( \pi\ga
\right)\Lambda(\ga)$ for some $1$-periodic trigonometric
polynomial $\Lambda$.
\end{itemize}
\el

\bp (a) : This follows by Lemma \ref{m-1}.

\noindent (b) : Note that $m_0(0)=\widetilde m_0(0)=1$ by the
assumptions in the general setup. This together with (a) and
\eqref{m-3} imples $m_0(1/2)=\widetilde m_0(1/2)=0$.

\noindent (c) : By (a) and Lemma \ref{m-4}, the functions $m_\ell,
\widetilde m_\ell$, $\ell=1,2\cdots, n$,
 can be factorized as
$$m_\ell(\ga)=e^{-\pi i \ga} \sin(\pi \ga) \Lambda_\ell(\ga), \
\widetilde m_\ell(\ga)=e^{-\pi i \ga} \sin(\pi \ga) \widetilde
\Lambda_\ell(\ga),$$ for some $1$-periodic trigonometric
polynomials $\Lambda_\ell, \widetilde \Lambda_\ell$. Combining
this with \eqref{m-2} leads to
$ 1-\overline{m_0(\ga)}\widetilde m_0(\ga)= \sin^2\left( \pi\ga
\right)\Lambda(\ga),$ where $ \Lambda(\ga)=\sum_{\ell=1}^n
\overline{\Lambda_\ell(\ga)}\widetilde \Lambda_\ell(\ga).$ \ep

In Section \ref{1812c} we will show that, on the other hand, the
assumptions (a), (b) and (c) in Lemma \ref{m-5} implies that the
wavelet systems $\{D^jT_k \psi_1\}_{j,k\in \mz}, \ \{D^jT_k
\widetilde{\psi_1}\}_{j,k\in \mz}$ can be extended to pairs of
dual wavelet frames by adding two wavelet systems.
\section{Extension with one pair of generators} \label{1812b}

In the rest of the paper we will consider scaling functions
$\varphi, \widetilde{\varphi}\in \ltr$  as in the general setup,
with associated trigonometric polynomial masks $m_0,
\widetilde{m_0}.$  Assuming that we have given trigonometric
polynomials $m_1, \widetilde{m_1}$ and defined the associated
functions $\psi_1, \widetilde{\psi_1}$ by \eqref{r-1}, our goal is
to extend the Bessel sequences $\{D^jT_k \psi_1\}_{j,k\in \mz}, \
\{D^jT_k \widetilde{\psi_1}\}_{j,k\in \mz}$ to pairs of dual
wavelet frames. It turns out to be convenient to consider the
functions $\widetilde M_{\alpha}$ and $\widetilde M_{\beta},$
defined by
\begin{eqnarray} \label{1812d}
\widetilde M_{\alpha}(\ga)&:=&1-\overline{m_0(\ga)} {\widetilde
m}_0(\ga)- \overline{m_1(\ga)} {\widetilde m}_1(\ga); \\
\label{1812e} \widetilde M_{\beta}(\ga)&:=&- \overline{m_0(\ga)}
{\widetilde m}_0(\ga+1/2)- \overline{m_1(\ga)} {\widetilde
m}_1(\ga+1/2).
\end{eqnarray}

\bl If the conditions (a), (b) and (c) in Lemma \ref{m-5} are
satisfied, then $\widetilde M_\alpha$ and $\widetilde M_\beta$ can
be factorized as
\begin{equation} \label{m-6}
   \widetilde M_\alpha(\ga)=\sin^2(\pi \ga) \Lambda_\alpha(\ga), \  \ \
\widetilde M_\beta(\ga)=-i \sin(\pi \ga)\cos(\pi \ga) \Lambda_\beta(\ga)
\end{equation}
for some $1$-periodic trigonometric polynomials $\Lambda_\alpha$
and $\Lambda_\beta$. \el

\bp By Lemma \ref{m-4}, and (a) and (b) of Lemma \ref{m-5},
$m_\ell, \widetilde m_\ell$, $\ell=0,1$ can be factorized as
\begin{eqnarray*}
m_0(\ga)=e^{-\pi i \ga} \cos(\pi \ga) \Lambda_0(\ga), &&
\widetilde m_0(\ga)=e^{-\pi i \ga} \cos(\pi \ga) \widetilde \Lambda_0(\ga),\\
m_1(\ga)=e^{-\pi i \ga} \sin(\pi \ga) \Lambda_1(\ga), &&
\widetilde m_1(\ga)=e^{-\pi i \ga} \sin(\pi \ga) \widetilde
\Lambda_1(\ga),
\end{eqnarray*}
for some $1$-periodic trigonometric polynomials $\Lambda_\ell,
\widetilde \Lambda_\ell$ for $\ell=0,1.$ Together with (c) of
Lemma \ref{m-5} this implies
\begin{eqnarray*}
   \widetilde M_\alpha(\ga) = \sin^2(\pi \ga)
   \left( \Lambda(\ga)-\overline{\Lambda_1(\ga)}\widetilde \Lambda_1(\ga)  \right)
   =: \sin^2(\pi \ga)\Lambda_\alpha(\ga),
\end{eqnarray*}
and
\begin{eqnarray*}
  \widetilde M_\beta(\ga)&=& - i\sin(\pi \ga)\cos(\pi \ga)
  \left(\overline{\Lambda_0(\ga)}\widetilde \Lambda_0(\ga+\pi/2)-
\overline{\Lambda_1(\ga)}\widetilde \Lambda_1(\ga+\pi/2)
   \right)\\
   &=:& - i \sin(\pi \ga)\cos(\pi \ga) \Lambda_\beta(\ga),
\end{eqnarray*}
as desired.\ep

We are now ready to state a condition for extension of  MRA-type
wavelet systems $\{D^jT_k\psi_1\}_{j, k \in \mz},
\{D^jT_k\widetilde \psi_1\}_{j, k \in \mz}$ to dual frames
$\{D^jT_k\psi_\ell\}_{\ell=1,2;j, k \in \mz}, \
\{D^jT_k\widetilde \psi_\ell\}_{\ell=1,2;j, k \in \mz}.$ Note
that in Theorem \ref{2012-14} below, the condition (i) in (I)
means that $\{D^jT_k\psi_2\}_{j, k \in \mz}$ and
$\{D^jT_k\widetilde \psi_2\}_{j, k \in \mz}$ are Bessel sequences,
while (ii) simply is the MEP-condition. We also note that the
proof shows how to choose the corresponding masks $m_2,
\widetilde{m_2}.$ With this information we can find the functions
$\psi_2, \widetilde{\psi_2}$ explicitly: if, e.g., $m_2(\ga)= \sum
d_ke^{2\pi i k\ga},$ then $\psi_2= 2\sum d_k \varphi(2x+k).$

\begin{thm}\label{2012-14}
Let $\varphi, \widetilde \varphi \in L^2(\mr)$ be  as in the
general setup, with trigonometric polynomial masks $m_0,
\widetilde m_0,$ respectively. Let $m_1, \widetilde m_1 $ be
trigonometric polynomials, and define $\psi_1, \widetilde
\psi_1\in L^2(\mr)$ by \eqref{r-1}. Assume that the conditions
(a), (b) and (c) in Lemma \ref{m-5} for $m_0, \widetilde m_0, m_1,
\widetilde m_1$ are satisfied. Then the following are equivalent:
\begin{itemize}
    \item[{\rm (I)}]
There exist $1$-periodic trigonometric polynomials $m_2,
\widetilde m_2$ such that
\begin{itemize}
    \item[{\rm (i)}] $m_2(0)=\widetilde m_2(0)=0$;
    \item[{\rm (ii)}] the matrix-valued functions $M, \widetilde M$ in
\eqref{m-12} with $n=2$ satisfy
\begin{equation*}
\widetilde M(\ga)^*M(\ga)=I,\  \ga \in \mt,
\end{equation*}
\end{itemize}
    \item[{\rm (II)}]
$\widetilde M_{\alpha}(\ga)\widetilde M_{\alpha}(\ga+1/2)
=\widetilde M_{\beta}(\ga)\widetilde M_{\beta}(\ga+1/2),\ \ga\in
\mt. $
\end{itemize}
In the affirmative case, the multi-wavelet systems
$\{D^jT_k\psi_\ell\}_{\ell=1,2; j, k \in \mz}$ and
$\{D^jT_k\widetilde \psi_\ell\}_{\ell=1,2; j, k \in \mz}$, with
$\psi_2, \widetilde \psi_2$ defined by \eqref{r-1}, form dual
frames  for $L^2(\mr)$.
\end{thm}
\bp (I)$\Rightarrow$(II): This follows from \eqref{m-2} and
\eqref{m-3}:
\begin{eqnarray*}
\widetilde M_{\alpha}(\ga)\widetilde M_{\alpha}(\ga+1/2)
&=& \overline{m_2(\ga)}\widetilde m_2(\ga)\overline{m_2(\ga+1/2)}{\widetilde m}_2(\ga +1/2)\\
&=& \widetilde M_{\beta}(\ga)\widetilde M_{\beta}(\ga+1/2).
\end{eqnarray*}

\noindent (II)$\Rightarrow$(I): Let $\widetilde M_\alpha,
\widetilde M_\beta$ be factorized as \eqref{m-6} for $1$-periodic
trigonometric polynomials $\Lambda_\alpha, \Lambda_\beta$. Then
Lemma \ref{m-5} (b) implies
\begin{equation}\label{m-8}
  \Lambda_\alpha(\ga)\Lambda_\alpha(\ga+1/2)=\Lambda_\beta(\ga)\Lambda_\beta(\ga+1/2).
\end{equation}
Let $\Gamma$ be the common factor of $\Lambda_\alpha$ and
$\Lambda_\beta$, that is,
\begin{equation} \label{m-10}
\Lambda_\alpha(\ga)= \Gamma(\ga) \Gamma_\alpha(\ga), \
\Lambda_\beta(\ga)= \Gamma(\ga) \Gamma_\beta(\ga),
\end{equation}
for some $1$-periodic trigonometric polynomials $\Gamma_\alpha,
\Gamma_\beta$ with no common factors. This together with
\eqref{m-8} implies
$$\Gamma_\alpha(\ga)\Gamma_\alpha(\ga+1/2)=\Gamma_\beta(\ga)\Gamma_\beta(\ga+1/2).$$
Since $\Gamma_\alpha$ and $\Gamma_\beta$ have no common factor, we
have
\begin{equation} \label{m-9}
\Gamma_\beta(\ga)=\Gamma_\alpha(\ga +1/2)
\end{equation}
up to constant. Define $m_2$ and $\widetilde m_2$ by
\begin{equation*}
   m_2(\ga):= e^{-\pi i\ga} \sin (\pi \ga) \overline{\Gamma(\ga)},\
   \widetilde m_2(\ga):= e^{-\pi i\ga} \sin (\pi \ga) \Gamma_\alpha(\ga).
\end{equation*}
Then (i) is trivial. By \eqref{m-10} and \eqref{m-9}, we have
$$\overline{m_2(\ga)}\widetilde m_2(\ga)= \sin^2(\pi \ga)\Gamma(\ga)\Gamma_\alpha(\ga)
= \widetilde M_\alpha(\ga),
$$
and
\begin{eqnarray*}
\overline{m_2(\ga)}\widetilde m_2(\ga+1/2)=
 -i \sin(\pi \ga) \cos(\pi \ga) \Gamma(\ga) \Gamma_\alpha(\ga+1/2)
= \widetilde M_\beta(\ga).
\end{eqnarray*}
These lead to (ii) by \eqref{m-2} and \eqref{m-3} with $n=2$.\ep

For the case where $\varphi= \widetilde{\varphi}=B_2$, we can
characterize the possible trigonometric polynomials $m_1,
\widetilde{m_1}$ with at most three terms that satisfy the
conditions in Theorem \ref{2012-14}. Our main reason for stating
this is that we can use the result to identify concrete candidates
for pairs of wavelet systems that can not be extended to a pair of
dual wavelet frames by adding a single pair of wavelet systems,
see Example \ref{2612a}.

\bc\label{m-25} Let $d_0,d_1, \widetilde d_0, \widetilde d_1 \in
\mc.$ Define $\psi_1$ and $\widetilde \psi_1$ by
\begin{eqnarray}
   \psi_1(x)&:=&d_0 B_2(2x)+(d_1-d_0)B_2(2x-1)-d_1 B_2(2x-2); \label{m-26}\\
   \widetilde \psi_1(x)&:=&\widetilde d_0 B_2(2x)+(\widetilde d_1-\widetilde d_0)B_2(2x-1)-\widetilde d_1
   B_2(2x-2).\label{m-27}
\end{eqnarray}
Then the following are equivalent:
\begin{itemize}
    \item[{\rm(a)}] There exist $1$-periodic trigonometric polynomials $m_2$ and $\widetilde m_2$ such that
    (I) in Theorem \ref{2012-14} holds;
    \item[{\rm(b)}] $3\overline d_0 \widetilde d_0 +3 \overline d_1 \widetilde d_1  - \overline d_1 \widetilde d_0  -\overline d_0 \widetilde d_1 =2$.
\end{itemize}

\ec

\bp  Let $\varphi:=\widetilde \varphi:=B_2$ with the associated
masks $m_0(\ga):=\widetilde m_0(\ga):=\left( e^{-\pi i \ga} \cos(\pi \ga) \right)^2.$
From \eqref{m-26} and \eqref{m-27}, $m_1$ and $\widetilde m_1$ are
defined by
\begin{equation}\label{m-30}
m_1(\ga)= \frac{d_0}{2}+\frac{d_1- d_0}{2}e^{-2\pi i \ga}
-\frac{d_1}{2}e^{-4\pi i \ga}, \ \  \widetilde m_1(\ga)=
\frac{\widetilde d_0}{2}+\frac{\widetilde d_1- \widetilde
d_0}{2}e^{-2\pi i \ga} -\frac{\widetilde d_1}{2}e^{-4\pi i \ga}.
\end{equation}
Trivially,  (a), (b) and (c) in Lemma \ref{m-5} for $m_0,
\widetilde m_0, m_1, \widetilde m_1$ are satisfied. We now check
the condition (II) in Theorem \ref{2012-14}. Note that
$$
\widetilde M_\alpha(\ga)= 1-\cos^4(\pi \ga) -\overline{m_1(\ga)} \widetilde m_1(\ga),\ \
\widetilde M_\beta(\ga)= \cos^2(\pi \ga)\sin^2(\pi \ga) -
\overline{m_1(\ga)} \widetilde m_1(\ga+1/2).
$$
Then we have
\begin{eqnarray*}
   \widetilde M_\alpha(\ga)\widetilde M_\alpha(\ga+1/2)
   &=&\left(1-\cos^4(\pi \ga) \right) \left(1-\sin^4(\pi \ga) \right)
    -  \left(1-\cos^4(\pi \ga) \right)\overline{m_1(\ga+1/2)} \widetilde m_1(\ga+1/2)\\
&& -  \left(1-\sin^4(\pi \ga) \right)\overline{m_1(\ga)} \widetilde m_1(\ga)
+\overline{m_1(\ga)} \overline{m_1(\ga+1/2)} \widetilde
m_1(\ga) \widetilde m_1(\ga+1/2)
\end{eqnarray*}
and
\begin{eqnarray*}
   && \widetilde M_\beta(\ga)\widetilde M_\beta(\ga+1/2) \\
   &&= \cos^4(\pi \ga) \sin^4(\pi \ga)
    - \sin^2(\pi \ga) \cos^2(\pi \ga)\overline{m_1(\ga+1/2)} \widetilde m_1(\ga)\\
   && - \cos^2(\pi \ga) \sin^2(\pi \ga)\overline{m_1(\ga)} \widetilde m_1(\ga+1/2)
   +\overline{m_1(\ga)} \overline{m_1(\ga+1/2)} \widetilde m_1(\ga) \widetilde m_1(\ga+1/2).
\end{eqnarray*}
Using the identities
\begin{eqnarray*}
  &&1-\cos^4(\pi \ga)=\sin^2(\pi \ga) \left(1+\cos^2(\pi \ga) \right), \ \
1-\sin^4(\pi \ga)=\cos^2(\pi \ga) \left(1+\sin^2(\pi \ga) \right),\\
  &&  \cos^2(\pi \ga)-\cos^4(\pi \ga)=\sin^2(\pi \ga)-\sin^4(\pi \ga)
  =\cos^2(\pi\ga)\sin^2(\pi\ga),
\end{eqnarray*}
the condition (II) in Theorem \ref{2012-14}  is equivalent to
\begin{eqnarray}
&&  \sin^2(\pi\ga) \cos^2(\pi \ga)
  \left(2- \left(\overline{m_1(\ga)}-\overline{m_1(\ga+1/2)} \right)
  \left( \widetilde m_1(\ga)- \widetilde m_1(\ga+1/2) \right) \right) \nonumber \\
&&=\sin^2(\pi \ga)\overline{m_1(\ga+1/2)}\widetilde m_1(\ga+1/2)+
\cos^2(\pi \ga)\overline{m_1(\ga)}\widetilde m_1(\ga)
\label{m-23}.
\end{eqnarray}
From \eqref{m-30}, $m_1$ and $\widetilde m_1$ can be factorized as
\begin{eqnarray*}
m_1(\ga)=e^{-\pi i\ga}\sin(\pi \ga) \Lambda_1(\ga), \ \
\widetilde m_1(\ga)=e^{-\pi i\ga}\sin(\pi \ga) \widetilde
\Lambda_1(\ga),
\end{eqnarray*}
where
\begin{equation} \label{m-31}
\Lambda_1(\ga):=i(d_0+ d_{1} e^{-2\pi i\ga} ), \ \
\widetilde \Lambda_1(\ga):=i( \widetilde d_0+\widetilde d_{1}
e^{-2\pi i\ga}).
\end{equation}
Since
$\overline{m_1(\ga)}\widetilde m_1(\ga)=\sin^2(\pi \ga)\overline{\Lambda_1(\ga)} \widetilde \Lambda_1(\ga), $
 \eqref{m-23} is equivalent with
\begin{eqnarray*}
 2- \left(\overline{m_1(\ga)}-\overline{m_1(\ga+1/2)} \right)
  \left( \widetilde m_1(\ga)- \widetilde m_1(\ga+1/2) \right) =\overline{\Lambda_1(\ga+1/2)} \widetilde \Lambda_1(\ga+1/2)
+\overline{\Lambda_1(\ga)} \widetilde \Lambda_1(\ga),
\end{eqnarray*}
This together with \eqref{m-30} and \eqref{m-31} leads to
$2- \left(\overline d_1 -\overline d_0\right)
\left( \widetilde d_1 -\widetilde d_0\right)
= 2 (\overline d_0 \widetilde d_0 + \overline d_1 \widetilde d_1),$
that is, $ 3 \overline d_0 \widetilde d_0 +3 \overline d_1
\widetilde d_1- \overline d_1 \widetilde d_0-\overline d_0
\widetilde d_1=2,$ as desired.  \ep

Let us consider a concrete case and find the functions $\psi_2,
\widetilde{\psi_2}$ explicitly:
\begin{ex}
Let $d_0=1, d_1=0, \widetilde d_0=\frac{1}{2}, \widetilde
d_1=-\frac{1}{2}$ and consider
\begin{eqnarray*}
  \psi_1 (x) &=& B_2(2x)-B_2(2x-1),
  \widetilde \psi_1(x) = \frac{1}{2} B_2(2x)- B_2(2x-1) +\frac{1}{2} B_2(2x-1),
  \end{eqnarray*}
with the associated masks
\begin{eqnarray*}
  m_1(\ga)= i e^{-\pi i \ga}\sin(\pi \ga),\ \
  \widetilde m_1(\ga)= - e^{-2\pi i \ga}\sin^2(\pi \ga).
\end{eqnarray*}
Then (b) in Corollary \ref{m-25} holds. We now use
(II)$\Rightarrow  $(I) in Theorem \ref{2012-14} to construct
$\psi_2$ and $\tilde \psi_2$ such that the multi-wavelet systems
$\{D^jT_k\psi_\ell\}_{\ell=1,2; j, k \in \mz}$ and
$\{D^jT_k\widetilde \psi_\ell\}_{\ell=1,2; j, k \in \mz}$ form
dual frames  for $L^2(\mr)$. A direct calculation shows that
\begin{eqnarray*}
  \widetilde M_\alpha(\ga)= \sin^2(\pi \ga) \Lambda_\alpha(\ga), \ \
  \widetilde M_\beta(\ga)= -i \sin(\pi \ga)\cos(\pi \ga) \Lambda_\beta(\ga),
\end{eqnarray*}
where
\begin{eqnarray*}
\Lambda_\alpha(\ga)= \left(1+e^{-2\pi i \ga} \right)
  \left(e^{2\pi i \ga} +3  \right), \ \
\Lambda_\beta(\ga)= \left(1+e^{-2\pi i \ga} \right)
  \left(e^{2\pi i \ga} -3\right)
\end{eqnarray*}
Letting $\Gamma(\ga) = -i \left(1+e^{-2\pi i \ga} \right)$ we have
$\Lambda_\alpha(\ga)= \Gamma(\ga) \Gamma_\alpha(\ga),
\Lambda_\beta(\ga)= \Gamma(\ga) \Gamma_\beta(\ga),$
where $\Gamma_\alpha(\ga)= i \left(e^{2\pi i \ga}+3 \right),
\Gamma_\beta(\ga)= i \left(e^{2\pi i \ga}-3 \right).$
Let
\begin{eqnarray*}
&& m_2(\ga)= e^{-\pi i \ga}\sin(\pi \ga)\overline{\Gamma(\ga)}= e^{-\pi i \ga}\sin(\pi \ga) i \left( e^{2\pi i \ga}+1 \right),  \\
&&  \widetilde m_2(\ga)= e^{-\pi i \ga}\sin(\pi \ga) \Gamma_\alpha(\ga)=
   e^{-\pi i \ga}\sin(\pi \ga) i \left( e^{2\pi i \ga}+3 \right).
\end{eqnarray*}
and define $\psi_2, \widetilde \psi_2$ by \eqref{r-1}, i.e.,
\begin{eqnarray*}
 \psi_2(x) = B_2(2x+1)-B_2(2x-1), \ \
   \widetilde \psi_2(x)=  B_2(2x+1)+2 B_2(2x) -3 B_2(2x-1).
  \end{eqnarray*}
By Theorem \ref{2012-14}, $\{D^jT_k\psi_\ell\}_{\ell=1,2; j, k \in
\mz}$ and $\{D^jT_k\widetilde \psi_\ell\}_{\ell=1,2; j, k \in
\mz}$ form dual frames.                  \ep
\end{ex}

Remember that our main motivation is the question whether any pair
of Bessel wavelet systems can be extended to a pair of dual
wavelet frames by adding a single pair of wavelet systems. Based
on Corollary \ref{m-25} we can state a concrete candidate for a
counter example:

\begin{ex} \label{2612a}
 Let $d_0=\widetilde d_0=1, d_1=\widetilde d_1=0$ and consider
$$\psi_1(x)=\widetilde \psi_1(x)=B_2(2x)-B_2(2x-1).$$
Then we have
$$3\overline d_0 \widetilde d_0 +3 \overline d_1 \widetilde d_1  - \overline d_1 \widetilde d_0  -\overline d_0 \widetilde d_1 =3\neq 2.$$
By Corollary \ref{m-25}, there does not exist a pair of
$\psi_2,\widetilde \psi_2$ such that $\widetilde M(\ga)^*M(\ga)=I$
with $n=2$. Thus the wavelet systems $\{D^jT_k\psi_1\}_{j, k \in
\mz}$ and $\{D^jT_k\widetilde \psi_1\}_{j, k \in \mz}$ can not be
extended to dual pairs of MRA-based wavelet frames
$\{D^jT_k\psi_\ell\}_{\ell=1,2; j, k \in \mz}, \
\{D^jT_k\widetilde \psi_\ell\}_{\ell=1,2; j, k \in \mz}.$ On the
other hand it is an open question whether the systems
$\{D^jT_k\psi_1\}_{j, k \in \mz}, \ \{D^jT_k\widetilde
\psi_1\}_{j, k \in \mz}$ can be extended to  dual wavelet pairs
$\{D^jT_k\psi_\ell\}_{\ell=1,2; j, k \in \mz}, \
\{D^jT_k\widetilde \psi_\ell\}_{\ell=1,2; j, k \in \mz}$
for which the functions $\psi_2, \widetilde{\psi_2}$ do not have the MRA-structure. \ep
\end{ex}

The following example illustrates the flexibility of the
MEP-approach compared with the UEP. In fact, we will consider
masks $m_0$ and $m_1$ as in the general setup, with the property
that the wavelet system $\{D^jT_k\psi_1\}_{ j, k \in \mz}$ can not
be extended to a tight wavelet frame
$\{D^jT_k\psi_\ell\}_{\ell=1,2; j, k \in \mz}$ using the UEP. On
the other hand, by considering a second pair of masks
$\widetilde{m_0}, \widetilde{m_1},$ the wavelet systems
$\{D^jT_k\psi_1\}_{ j, k \in \mz}$ and
$\{D^jT_k\widetilde{\psi_1}\}_{ j, k \in \mz}$ can be extended to
dual pairs $\{D^jT_k\psi_\ell\}_{\ell=1,2; j, k \in \mz}, \
\{D^jT_k\widetilde{\psi_\ell}\}_{\ell=1,2; j, k \in \mz}$ via the
MEP:

\begin{ex}\label{m-21}
Consider
$$ \varphi(x):=B_1(x), \ \
\psi_1(x):=\frac{1}{4}B_1(2x+2)+\frac{1}{4}B_1(2x+1)
-\frac{1}{4}B_1(2x)-\frac{1}{4}B_1(2x+2)$$ with the associated
masks
$$m_0(\ga)= e^{-\pi i  \ga} \cos(\pi \ga), \ \
m_1(\ga)=i e^{\pi i \ga } \sin(\pi \ga)\cos^2(\pi \ga).$$ Let
\begin{eqnarray*}
M_{\alpha}(\ga)&=&1-\left|m_0(\ga)\right|^2 -\left|m_1(\ga)\right|^2, \ \
M_{\beta}(\ga)=- m_0(\ga)\overline{m}_0(\ga+1/2)-
m_1(\ga)\overline{m}_1(\ga+1/2).
\end{eqnarray*}
A direct calculation shows that
$ M_{\alpha}(\ga) M_{\alpha}(\ga+1/2)\not\equiv M_{\beta}(\ga) M_{\beta}(\ga+1/2);$
by \cite[Theorem 2.2]{CKK-CMB}, this implies that there does not exist $m_2\in
L^\infty(\mr)$ such that the multi-wavelet system
$\{D^jT_k\psi_\ell\}_{\ell=1,2; j, k \in \mz}$, with $\psi_2$
defined by \eqref{r-1}, forms a Parseval frame  for $L^2(\mr)$.

On the other hand, consider now
$$\widetilde \varphi(x):=B_3(x+1), \ \
\widetilde \psi_1 (x):=B_3(2x+2) -B_3(2x+1)$$ with the associated
masks
$  \widetilde m_0(\ga) = e^{-\pi i \ga} \cos^3(\pi \ga),
\widetilde m_1(\ga) = i e^{\pi i \ga} \sin(\pi \ga).$ A direct
calculation shows that
$$\widetilde M_{\alpha}(\ga)=\sin^2(\pi\ga), \ \
\widetilde M_{\beta}(\ga)=-i\sin(\pi \ga)\cos(\pi \ga).$$
Then we
have
$$\widetilde M_{\alpha}(\ga)\widetilde M_{\alpha}(\ga+1/2)
=\widetilde M_{\beta}(\ga)\widetilde M_{\beta}(\ga+1/2)=\sin^2(\pi
\ga)\cos^2(\pi \ga). $$ It is easy to check that the conditions
(a), (b) and (c) in Lemma \ref{m-5} are satisfied, so by Theorem \ref{2012-14}
there exists functions $\psi_2, \widetilde{\psi_2}$ such that
$\{D^jT_k\psi_\ell\}_{\ell=1,2; j, k \in \mz}$ and
$\{D^jT_k\widetilde \psi_\ell\}_{\ell=1,2; j, k \in \mz}$ form
dual  frames. By the proof of Theorem \ref{2012-14} we can
define $\psi_2$ and $\widetilde \psi_2$ by
\begin{equation*}
   \psi_2(x)=B_1(2x)-B_1(2x-1), \ \
   \widetilde \psi_2(x)=B_3(2x+1)-B_3(2x)
\end{equation*}
with the associated masks
$m_2(\ga)=\widetilde m_2(\ga)= i e^{-\pi i\ga} \sin (\pi \ga).$\ep

\end{ex}

\section{Extension with two pairs of generators} \label{1812c}

For the sake of completeness we now prove that if the necessary
conditions in Lemma \ref{m-5} are satisfied, then we can always
extend $\{D^jT_k\psi_1  \}_{j, k \in \mz}$,
$\{D^jT_k\widetilde\psi_1\}_{j, k \in \mz}$ to dual wavelet frame
$\{D^jT_k\psi_\ell  \}_{\ell=1,2,3;j, k \in \mz}$,
$\{D^jT_k\widetilde\psi_\ell\}_{\ell=1,2,3;j, k \in \mz},$ by
adding two pairs of wavelet systems with MRA-structure. Note that
if it is not required that the functions $\psi_2, \psi_3,
\widetilde{\psi_2}, \widetilde{\psi_3}$ have MRA-structure, it is
trivial that such an extension exists: we can just take
$\psi_2=\psi_1, \widetilde{\psi_2}=-\widetilde{\psi_1}$ to
``cancel" the action of the given wavelet system, and then take
$\psi_3$  and $\widetilde{\psi_3}$ to be any pair of functions
generating dual wavelet frames. But this argument does not work
under the assumptions in the current paper. For example, if we
consider the MRA-setup where the functions $\psi_\ell,
\widetilde{\psi_\ell}$ are generated by a B-spline, $\varphi=B_N,
N>1,$ then there does not exist dual wavelet frames
$\{D^jT_k\psi_3  \}_{j, k \in \mz}$,
$\{D^jT_k\widetilde\psi_3\}_{j, k \in \mz}$ where \bes \psi_3=
\sum c_k B_N(2x-k), \ \ \widetilde{\psi_3}= \sum \widetilde{c_k}
B_N(2x-k),\ens with finite coefficient sequences, i.e., the masks
$m_3, \widetilde{m_3}$ can not be trigonometric polynomials; see
\cite{CHS} (the result is repeated in Theorem 14.5.1 in \cite{CB}).
Thus, we have to provide a different argument:

\bt \label{m-19} Let $\varphi, \widetilde \varphi\in L^2(\mr)$ be
as in the general setup, with trigonometric polynomial masks
$m_0,\widetilde m_0.$ Let $m_1,\widetilde m_1$ be trigonometric
polynomials, and define $\psi_1, \widetilde \psi_1 $ by
\eqref{r-1}. Assume that the conditions (a), (b) and (c) in Lemma
\ref{m-5} for $m_0, \widetilde m_0, m_1, \widetilde m_1$ are
satisfied. Then  there exist  trigonometric polynomials $m_2,
m_3,\widetilde m_2, \widetilde m_3$ such that
 $\{D^jT_k\psi_\ell\}_{\ell=1,2,3; j, k \in \mz}$ and $\{D^jT_k\widetilde \psi_\ell\}_{\ell=1,2,3; j, k \in \mz}$,
with $\psi_2,\widetilde\psi_2, \psi_3, \widetilde\psi_3 $ defined
by \eqref{r-1}, form dual  frames  for $L^2(\mr)$.\et \bp We
construct $m_2, \widetilde m_2, m_3, \widetilde m_3$ such that
\begin{itemize}
    \item[(1)] $m_2(0)=\widetilde m_2(0)=m_3(0)=\widetilde m_3(0)=0$;
    \item[(2)] the matrix-valued functions $M,\widetilde M$ in \eqref{m-12} with $n=3$ satisfy
    $$\widetilde M(\ga)^*M(\ga)=I,\  \ga \in \mt.$$
\end{itemize} Choose $M_\alpha, M_\beta$ as in \eqref{1812d}--\eqref{1812e}.
We first define $m_2$ and $\widetilde m_2$ by
\begin{eqnarray}
m_2(\ga)&:=&\overline{\widetilde M_\alpha(\ga)}+\overline{\widetilde M_\beta(\ga)}
= \sin(\pi \ga) \left( \sin(\pi \ga)
\overline{\Lambda_\alpha(\ga)}
+i \cos(\pi \ga) \overline{\Lambda_\beta(\ga)}  \right) \label{m-14} \\
\widetilde m_2(\ga)&:=&\sin^2(\pi \ga), \label{m-15}
\end{eqnarray}
where $\Lambda_\alpha$ and $\Lambda_\beta$ are defined as in
\eqref{m-6}. Then
\begin{eqnarray*}
1-\sum_{\ell=0}^2 \overline{m_\ell(\ga)}\widetilde m_\ell(\ga)&=&
\widetilde M_\alpha(\ga)     - \overline{m_2(\ga)}\widetilde m_2(\ga)
= \widetilde M_\alpha(\ga)- \left(\widetilde M_\alpha(\ga)+\widetilde M_\beta(\ga)\right) \sin^2(\pi\ga) \\
&=& \sin^2(\pi \ga) \cos(\pi \ga) \left(  \cos(\pi \ga)
\Lambda_\alpha(\ga)+ i \sin(\pi \ga) \Lambda_\beta(\ga)
\right)
\end{eqnarray*}
and
\begin{eqnarray*}
 \sum_{\ell=0}^2 \overline{m_\ell(\ga)}\widetilde m_\ell(\ga+1/2)&=&
-\widetilde M_\beta(\ga) +  \overline{m_2(\ga)}\widetilde m_2(\ga+1/2)\\
&=&\sin^2(\pi \ga) \cos(\pi \ga) \left(  \cos(\pi \ga)
\Lambda_\alpha(\ga)+ i \sin(\pi \ga) \Lambda_\beta(\ga)
\right)
\end{eqnarray*}
We now define $m_3, \widetilde m_3$ by
\begin{eqnarray}
m_3(\ga)&:=&\sin(\pi \ga) \left( \cos(\pi
\ga)\overline{\Lambda_\alpha(\ga)}-
i \sin(\pi \ga) \overline{\Lambda_\beta(\ga)}  \right) \label{m-16}\\
\widetilde m_3(\ga)&:=& \sin(\pi \ga) \cos(\pi \ga). \label{m-17}
\end{eqnarray}
Then it is easy to see that
\begin{eqnarray*}
   && 1- \sum_{\ell=0}^2 \overline{m_\ell(\ga)}\widetilde m_\ell(\ga)=\overline{m_3(\ga)}\widetilde m_3(\ga), \\
&& \sum_{\ell=0}^2 \overline{m_\ell(\ga)}\widetilde
m_\ell(\ga+1/2)= -\overline{m_3(\ga)}\widetilde m_3(\ga+1/2),
\end{eqnarray*}
which is equivalent to the condition (2). By construction, (1)
holds. \ep

Let us illustrate Theorem \ref{m-19} by returning to a
construction presented in Example 2.5 in \cite{CKK-CMB}. In that
example, masks $m_0, m_1$ as in the general setup were considered,
and it was shown how to extend the Bessel sequence
$\{D^jT_k\psi_1\}_{ j, k \in \mz}$ to a tight frame
$\{D^jT_k\psi_\ell\}_{\ell=1,2,3; j, k \in \mz}$ using the UEP.
However, the masks defining the functions $\psi_2, \psi_3$ are not
trigonometric polynomials, so the functions $\psi_2, \psi_3$ are
not compactly supported. Using the MEP we can now show that it is
possible to extend $\{D^jT_k\psi_1\}_{ j, k \in \mz}$ to a pair of
dual wavelet frames with compactly supported generators. The
example deals again with the B-splines, for convenience shifted to
be centered around $x=0:$

\begin{ex}
Let $\ell \geq 2.$ Consider
$$ \varphi(x):=B_{2\ell}(x + \ell), \ \
\psi_1(x):=\sum_{k=-\ell}^\ell \frac{1}{2^{2\ell-1}}{2\ell \choose
\ell+k} B_{2\ell}(2x-k+\ell)
$$
with the associated masks
$$m_0(\ga)=  \cos^{2\ell} \left(\pi \ga \right), \ \
 m_1(\ga)=   \sin^{2 \ell} (\pi \ga).$$
In \cite[Example 2.5]{CKK-CMB} it was shown that by letting
$$ \lambda_1(\ga):=1-\left(\cos^{2\ell}(\pi \ga)+\sin^{2\ell}(\pi \ga)\right)^2, \
 \lambda_2(\ga):=1-\left(\cos^{2\ell}(\pi \ga)-\sin^{2\ell}(\pi \ga)\right)^2,$$
the masks $m_2, m_3\in L^\infty(\mt)$ defined by
\begin{equation*}
 m_2(\ga):= \left\{
\begin{array}{ll}
\left( \sqrt{\lambda_1(\ga)}+\sqrt{\lambda_2(\ga)} \right)/2,& \ga \in [-1/2,0]+\mz\\
\left( \sqrt{\lambda_1(\ga-1/2)}-\sqrt{\lambda_2(\ga-1/2)}
\right)/2,& \ga \in [0, 1/2]+\mz
\end{array}
\right.
\end{equation*}
and
\begin{equation*}
m_3(\ga):= \left\{
\begin{array}{ll}
\left( \sqrt{\lambda_1(\ga)}-\sqrt{\lambda_2(\ga)}\right)/2,& \ga \in [-1/2,0]+\mz\\
\left( \sqrt{\lambda_1(\ga-1/2)}+\sqrt{\lambda_2(\ga-1/2)}
\right)/2,& \ga \in [0,1/2]+\mz
\end{array}
\right.  ,
\end{equation*}
satisfy the UEP conditions. We note that these masks are not
trigonometric polynomials.

Using the MEP we can now provide an alternative construction using
trigonometric polynomials. Consider now $ \widetilde
\varphi:=\varphi, \  \widetilde \psi_1:=\psi_1$ with the
associated masks
$$\widetilde m_0(\ga):= m_0(\ga)= \cos^{2\ell} \left(\pi \ga \right), \ \
 \widetilde m_1(\ga):= m_1(\ga)=  \sin^{2 \ell} (\pi \ga).$$
A direct calculation shows that
$$\widetilde M_\alpha(\ga)=\sin^2(\pi \ga)\Lambda_\alpha(\ga),\
\widetilde M_\beta(\ga)=-i\sin(\pi \ga) \cos(\pi \ga)
\Lambda_\beta(\ga),$$ where
\begin{eqnarray*}
   \Lambda_\alpha(\ga) &=& \left(\sum_{k=0}^{2\ell -1} \cos^{2k}(\pi \ga)\right)-\sin^{4\ell-2}(\pi \ga), \ \
\Lambda_\beta(\ga) = -2i \sin^{2\ell-1}(\pi
\ga)\cos^{2\ell-1}(\pi \ga).
\end{eqnarray*}
Define trigonometric polynomials $m_2, \widetilde m_2, m_3,
\widetilde m_3$ as in \eqref{m-14}-\eqref{m-17}:
\begin{eqnarray*}
m_2(\ga) &:=& \sin(\pi \ga) \left( \sin(\pi \ga)
\overline{\Lambda_\alpha(\ga)} +i \cos(\pi \ga)
\overline{\Lambda_\beta(\ga)}  \right), \ \
\widetilde m_2(\ga):=\sin^2(\pi \ga); \\
m_3(\ga)&:=&\sin(\pi \ga) \left( \cos(\pi
\ga)\overline{\Lambda_\alpha(\ga)}-
i \sin(\pi \ga) \overline{\Lambda_\beta(\ga)}  \right),\ \
\widetilde m_3(\ga):= \sin(\pi \ga) \cos(\pi \ga).
\end{eqnarray*}
Then we see that the masks satisfy the conditions in Theorem
\ref{m-19}. Hence we conclude that
$\{D^jT_k\psi_\ell\}_{\ell=1,2,3; j, k \in \mz}$ and
$\{D^jT_k\widetilde \psi_\ell\}_{\ell=1,2,3; j, k \in \mz}$, with
$\psi_2,\widetilde\psi_2, \psi_3, \widetilde\psi_3 $ defined by
\eqref{r-1}, form dual  frames  for $L^2(\mr).$ By construction
all the generators are compactly supported. \ep
\end{ex}

\end{document}